\documentclass[11pt,a4paper]{article}
\usepackage[utf8]{inputenc}
\usepackage[russian,english]{babel}
\usepackage{amsmath,amssymb,latexsym,amsthm,enumitem,graphicx,setspace}
\usepackage[numbers,square,sort&compress]{natbib}
\usepackage[dvipsnames]{xcolor} 
\usepackage{cleveref}

\newtheorem{theorem}{Theorem}
\newtheorem{proposition}{Proposition}
\theoremstyle{definition}
\newtheorem{remark}{Remark}

\newcounter{numitem}[subsection]
\newcommand*{\newitem}{\refstepcounter{numitem}\medskip\textbf{\thenumitem. }}

\renewcommand{\tilde}{\widetilde}
\newcommand*{\ol}{\overline}
\newcommand*{\eps}{\varepsilon}
\newcommand*{\ffi}{\varphi}
\newcommand*{\EE}{\mathsf E}
\newcommand*{\PP}{\mathsf P}
\newcommand*{\bbR}{\mathbb R}
\newcommand*{\cA}{\mathcal A}
\newcommand*{\cF}{\mathcal F}
\newcommand*{\cV}{\mathcal V}
\newcommand{\FF}{\mathbb{F}}
\newcommand{\E}{\mathsf E}
\newcommand{\I}{I}
\DeclareMathOperator{\Law}{Law}

\title{Sequential tracking of an unobservable two-state Markov process under\ Brownian noise}
\author{Alexey Muravlev\thanks{Steklov Mathematical Institute of the Russian
Academy of Sciences. 8 Gubkina St., Moscow 119991, Russia. E-mail:
almurav@mi-ras.ru}
\and
Mikhail Urusov\thanks{The University of Duisburg--Essen, Faculty of
Mathematics. 9 Thea-Leymann-Str., Essen 45127, Germany. E-mail: mikhail.urusov@uni-due.de}
\and
Mikhail Zhitlukhin\thanks{Steklov Mathematical Institute of the Russian
Academy of Sciences. 8 Gubkina St., Moscow 119991, Russia. E-mail:
mikhailzh@mi-ras.ru. Corresponding author.}}
\date{3 August 2019}

\begin{document}
\maketitle

\begin{abstract}
We consider an optimal control problem, where a Brownian motion with drift is
sequentially observed, and the sign of the drift coefficient changes at jump
times of a symmetric two-state Markov process. The Markov process itself is
not observable, and the problem consist in finding a $\{-1,1\}$-valued
process that tracks the unobservable process as close as possible. We present
an explicit construction of such a process.


\end{abstract}

\section{Introduction}
We consider a problem of sequential tracking of a symmetric two-state Markov
process, the values of which cannot be seen directly but are only observed
with noise. The presence of noise is modelled by that this process appears as
a local drift coefficient of an observable Brownian motion with drift. Under
a tracking procedure we understand a two-state process such that its value should be
equal to the value of the unobservable process as often as possible,
provided that there is a penalty for frequent switching of the value of the
tracking process.

This problem can be viewed as a multiple quickest changepoint detection
problem. Recall that in standard changepoint detection problems, the goal is
to detect a (single) change of some characteristic of an observable process,
e.g., the drift of a Brownian motion. Single changepoint problems are
well-studied in the literature and there exist many different settings and
methods of their solution (see, e.g., the monographs \cite{Shiryaev19,
PoorHadjiliadis08, TartakovskyNikiforov14}). Problems, where changes may
occur several times, are less investigated, though for a wide range of
applications multiple changepoint models seem to be more adequate; in
particular, it is interesting to note that when the theory of quickest
detection started to actively develop in the 1950-60s, multiple changepoint
settings were though to be the ``right ones'' for practical applications ---
see a historical review in \cite{Shiryaev10}.

We solve our problem by first reducing it to an optimal control problem for
the posterior mean of the unobservable process. Using standard filtering
techniques for Brownian motion, it is possible to write an explicit SDE for
the posterior mean process. This allows to obtain an optimal control
problem, where the control process assumes just two values, i.e., we get a
so-called optimal switching problem. Its solution is obtained by considering
a free-boundary problem for an ODE associated with the infinitesimal
operator of the posterior mean process. Although the solution of the latter
problem cannot be expressed in elementary functions, it is possible to
characterize it in a rather convenient form, which allows to find it
numerically.

Let us briefly mention other results in the literature related to our paper.
A similar multiple changepoint detection problem was studied by
\citet{Gapeev15}. His optimality criterion is somewhat different from ours,
and he considers general (non-symmetric) two-state Markov processes. The
main results of \cite{Gapeev15} consist in reduction of the changepoint
detection problem to coupled optimal stopping problems and, further, to
coupled free-boundary problems. Then the optimal boundaries at which
switching occurs are identified as unique functions satisfying the smooth
fit conditions. It seems to be difficult to find an explicit solution, but,
nevertheless, the paper establishes some analytic estimates for it.

\citet{BayraktarLudkovski09} considered a tracking problem for a compound
Poisson process with local arrival rate and jump distribution depending on
the state of an unobservable Markov chain. That problem was also reduced to
a coupled optimal stopping problem, though, due to the nature of Poisson
process, the method used to solve it is somewhat different from techniques
applied to continuous processes.

There are many results related to optimal switching problems for stochastic
processes (optimal control problems, where control processes assume values
from a finite set). An exposition of the topic can be found, for
example, in Chapter~5 of~\cite{Pham09}. Among various results, let us
mention the paper of \citet{BayraktarEgami10}, where two-state switching
problems were considered. General results of that paper show how a solution
of a switching problem can be characterized as a solution of coupled optimal
stopping problems. The paper also includes several examples of explicit
solutions, but all of them are related to the situation when a controlled
process is a diffusion on $(0,\infty)$ with $\infty$ being the natural
boundary, and $0$ a natural or absorbing boundary (on the contrary, in our
problem, we have a diffusion with finite inaccessible entrance boundaries).
An explicit solution to a two-state switching problem for a geometric
Brownian motion was obtained by \citet{LyVathPham07}.

Let us also mention the paper of \citet{CaiRosenbaumTankov17a} (and the
subsequent paper \cite{CaiRosenbaumTankov17b}), which studies a general
tracking problem, where an observer needs to adjust a controlled process to
keep it close to an observable It\^o process. The main results are related
to the asymptotic analysis of the cost function when the costs are small.
Besides general results, the paper includes examples of their applications
to particular processes. Also, examples of applied models related to optimal
switching problems can be found in mathematical finance. Among others, we
can mention the papers \cite{BrekkeOksendal94,DuckworthZervos01} on optimal
investment decisions; see also references in \cite{CaiRosenbaumTankov17a}.

Our paper is organized as follows. \Cref{sec:prob} contains a formal
statement of the problem. The main theorem describing its solution is stated
in~\Cref{sec:solution}, together with a discussion on how a numerical
procedure can be implemented. All the proofs are assembled
in~\Cref{sec:proofs}. Finally, the appendix contains analysis of the
boundary behavior of the posterior mean process and establishes its
auxiliary properties needed in the paper.

\section{The problem}\label{sec:prob}
We consider independent processes $W=(W_t)_{t\ge0}$ and
$\theta=(\theta_t)_{t\ge0}$ on a probability space $(\Omega,\cF,\PP_x)$,
where $W$ is a standard Brownian motion with $W_0=0$, and $\theta$ is a
c\`adl\`ag $\{-1,1\}$-valued continuous-time Markov process, which jumps
from $-1$ to $1$ and back with intensity $\lambda>0$, with $\EE_x
\theta_0=x\in[-1,1]$ (that is, under $\PP_x$ the random variable $\theta_0$
has the distribution $\frac{1+x}2 \delta_1+\frac{1-x}2 \delta_{-1}$).
Neither of the processes $W$ and $\theta$ is considered observable. Given
some $\mu>0$ we sequentially observe the process
\[
X_t=\int_0^t \mu\theta_s\,ds+W_t,\quad t\ge0,
\]
and the aim is to track the ``hidden signal'' $\theta$ as close as possible.
More precisely, we consider the right-continuous filtration
$\FF=(\cF_t)_{t\ge0}$,
$$
\cF_t=\bigcap_{\eps>0}\sigma(X_s;\,s\in[0,t+\eps]),\quad t\ge0,
$$
generated by $X$, define the space of controls $\cA$ as the space of
$\FF$-adapted $\{-1,1\}$-valued c\`adl\`ag processes that have finite number
of jumps on compact time intervals with a possibility to have a jump at time
zero (i.e., we distinguish between $A_{0-}$ and $A_0$, which will be
important in \eqref{eq:04122015a1.5}--\eqref{V} below), and, for some
$\alpha>0$ and $c_1,c_2\ge 0$ such that $c_1+c_2>0$, define the cost
function of the arguments $x\in[-1,1]$, $A\in\cA$,
\begin{multline}\label{J}
J(x,A)=\EE_x\biggl[
\int_0^\infty e^{-\alpha t} I(A_t\ne\theta_t)\,dt
\\+\sum_{t\ge0} e^{-\alpha t} I(\Delta A_t\ne0) (c_1 + c_2 \I(A_t\neq \theta_t))
\biggr]
\end{multline}
and the value function
\begin{equation}\label{eq:04122015a1}
\cV^*(x)=\inf_{A\in\cA} J(x,A).
\end{equation}
That is, we have the running cost of intensity $1$ when the control process
$A$ differs from the (unobservable) signal $\theta$ and the cost for
switching between the levels $\pm1$ in the control $A$, which consists of
the fixed cost $c_1$ and the additional cost $c_2$ for switching at an
incorrect time. All costs being discounted with rate~$\alpha$.

The goal of the problem that we consider is, given the numbers
$\lambda,\mu,\alpha>0$, $c_1,c_2\ge0$ with $c_1+c_2>0$ as inputs, to find
the infimum in~\eqref{eq:04122015a1} for all $x\in[-1,1]$ and the
corresponding optimal control.

Let us make three simple observations regarding the structure of the
solution that will be used in what follows.

(i) For $a\in\{-1,1\}$, define
\begin{equation}\label{eq:04122015a1.5}
\cA_a=\{A\in\cA:A_{0-}=a\}
\end{equation}
and
\begin{equation}\label{V}
V^*(x,a)=\inf_{A\in\cA_a} J(x,A),\quad x\in[-1,1].
\end{equation}
Obviously, the solutions of switching problems in~\eqref{V} both for $a=-1$
and $a=1$ yield the solution of~\eqref{eq:04122015a1}. Indeed, we clearly
have
$$
\cV^*(x)=V^*(x,1)\wedge V^*(x,-1),\quad x\in[-1,1],
$$
and if $A^{*,a}$ are optimal controls in problem~\eqref{V} ($a\in\{-1,1\}$),
then the optimal control in~\eqref{eq:04122015a1} is $A^* = A^{*,1}$ if
$V^*(x,1)<V^*(x,-1)$, and $A^* = A^{*,-1}$ otherwise.

(ii) Each control process $A\in\cA$ ($\equiv\cA_{-1}\cup\cA_1$) can be
identified with a sequence of stopping times $(\tau_n)_{n=0}^\infty$ of the
filtration $\FF$ such that
\begin{equation}\label{eq:09122015a0}
0=\tau_0\le\tau_1<\tau_2<\ldots,\quad
\tau_n\to\infty,
\end{equation}
and in the case $A\in\cA_a$, $a\in\{-1,1\}$,
\begin{equation}\label{eq:09122015a0.5}
\begin{aligned}
&A_{0-}=a,\\
&A_t=a\text{ for }t\in[\tau_{2n},\tau_{2n+1}),\\
&A_t=-a\text{ for }t\in[\tau_{2n+1},\tau_{2n+2}),\quad n\ge0.
\end{aligned}
\end{equation}
Notice that the first inequality in~\eqref{eq:09122015a0} is not strict,
which allows the possibility of jump at time zero.

(iii) By the symmetry of the setting it is enough to solve~\eqref{V} only
for $a=1$, as it is easy to check that
$$
V^*(x,-1)=V^*(-x,1),\quad x\in[-1,1].
$$
For the statement regarding optimal controls, see \Cref{rem:symmetry} below.

\section{Solution of the problem}
\label{sec:solution}

\subsection{Main result}
We begin with introducing several auxiliary objects that will be needed to
provide a solution of the problem.

A key role will be played by the posterior mean process of $\theta_t$, which
we define as the $\FF$-adapted process
\begin{equation}\label{eq:31072019a1}
M_t=\EE_x(\theta_t|\cF_t),\quad t\ge0.
\end{equation}
From general results of the filtration theory for diffusion processes (see,
e.g., Theorem~9.1 in~\cite{LiptserShiryaev:01}), one can deduce that under
$\PP_x$, $x\in[-1,1]$, the process $M$ satisfies the SDE
\begin{equation}\label{eq:04122015a3.5}
dM_t=-2\lambda M_t\,dt+\mu (1-M_t^2)\,d\ol W_t,
\quad M_0=x,
\end{equation}
where the innovation process $\ol W$ is an $(\cF_t,\PP_x)$-Brownian motion.
Moreover, the relation between $\ol W$ and $X$ is described by the formula
$\ol W_t=X_t-\int_0^t \mu M_s\,ds$, which gives a possibility to express $M$
through the observable process $X$ in a practical realisation of the optimal
tracking rule.

Notice that $M$ is a pathwise unique solution of~\eqref{eq:04122015a3.5}, as
the coefficients in~\eqref{eq:04122015a3.5} are locally Lipschitz on $\bbR$
(see Theorem~5.2.5 in~\cite{KaratzasShreve:91}).
Expression~\eqref{eq:31072019a1} implies that, for $x\in[-1,1]$, the
solution $M$ of~\eqref{eq:04122015a3.5} is $[-1,1]$-valued, but we can say
more. By computing the scale function of $M$ inside $(-1,1)$ we establish
that the boundaries $\pm1$ are inaccessible; in particular, $M$ is
$(-1,1)$-valued whenever $x\in(-1,1)$. A further computation entails that
$\pm1$ are entrance boundaries for~\eqref{eq:04122015a3.5}, i.e., a solution
$M$ to~\eqref{eq:04122015a3.5} can be started in $\pm1$, which then
immediately enters $(-1,1)$ and never leaves $(-1,1)$. For more detail to
these points, see the appendix.

Let us now introduce the differential operator $L$, which is associated with
$M$ and the discounting factor $\alpha$ from~\eqref{J}, and acts on
$C^2$-functions $f\colon(-1,1)\to\bbR$ by
\begin{equation}\label{eq:dif_op}
Lf(x) = \frac{\mu^2}{2} (1-x^2)^2\,\frac{d^2f}{dx^2}(x)
-2\lambda x\,\frac{df}{dx}(x)-\alpha f(x),
\quad x\in(-1,1).
\end{equation}
Consider the second-order linear ODE
\begin{equation}\label{ode}
LV(x) = -\frac12(1-x),\quad x\in(-1,1).
\end{equation}
We will relate the solution of the problem $V^*(x,1)$ to a solution of a
certain free boundary problem for~\eqref{ode}. A straightforward calculation
shows that a particular solution of \eqref{ode} is
\[
\tilde V(x)=\frac{1}{2\alpha} - \frac{x}{2(2\lambda+\alpha)},
\quad x\in[-1,1].
\]
Solutions of the corresponding homogeneous ODE cannot be expressed in
elementary functions, but the following \namecref{pr:phi-properties} states
their properties that will be need further. We prove it in
\Cref{sec:proof-aux} by deducing from the general theory of one-dimensional
diffusions.

\begin{proposition}
\label{pr:phi-properties}
Consider the homogeneous ODE that corresponds to~\eqref{ode},
\begin{equation}\label{eq:hom_ode}
Lf(x)=0,\quad x\in(-1,1).
\end{equation}
The following claims are true:
\begin{enumerate}[itemsep=0mm,label=(\alph*),leftmargin=*,topsep=\lineskip]
\item this ODE has a strictly decreasing strictly positive solution~$\ffi$;
\item every decreasing strictly positive solution to~\eqref{eq:hom_ode} is a
constant multiple of~$\ffi$;
\item $\ffi(-1):=\lim\limits_{x\to-1}\ffi(x)=+\infty$ and
$\ffi(1):=\lim\limits_{x\to1}\ffi(x)>0$.
\end{enumerate}
\end{proposition}

Given any such solution $\ffi$ we obtain a strictly increasing strictly
positive solution $\psi$ to~\eqref{eq:hom_ode} by setting
$\psi(x)=\ffi(-x)$. Clearly, $\psi$ and $\ffi$ are linearly independent,
hence we obtain a general solution to the inhomogeneous ODE~\eqref{ode} in
the form
\begin{equation}\label{eq:08122015a1}
V(x)=\tilde V(x)+K_1 \ffi(x)+K_2 \psi(x),\quad x\in(-1,1),
\end{equation}
where $K_1,K_2\in\bbR$. It follows from~\eqref{eq:08122015a1} together with
claim~(c) that the set of solutions $V$ to~\eqref{ode} which are bounded in
a left neighbourhood of~$1$~is
\begin{equation}\label{eq:09122015a5}
\{\tilde V-K\ffi:K\in\bbR\}.
\end{equation}

Now we are ready to formulate the main result.

\begin{theorem}\label{th:main}
(i) Suppose $c_1<(2\lambda + \alpha)^{-1}$. Then there exists a unique pair
\begin{equation}\label{eq:09122015a5.5}
(K,B)\in(0,\infty)\times\biggl(\frac{2c_1+c_2}{2(2\lambda+\alpha)^{-1}+c_2},\;
1\biggr)
\end{equation}
such that, with $V:=\tilde V-K\ffi$ (cf.~(\ref{eq:09122015a5})), we have
\begin{align}
\label{eq:09122015a6}
V(-B)&=V(B)+c_1 + \frac{c_2}{2}(1-B),\\
\label{eq:09122015a7}
V'(-B)&=-V'(B) + \frac{c_2}{2}.
\end{align}
Furthermore, the value function $V^*(x,1)$, $x\in[-1,1]$, is given by
\begin{equation}\label{eq:09122015a8}
V^*(x,1)=\begin{cases}
V(x),&x\in[-B,1],\\
V(-x)+c_1 + \dfrac{c_2}{2}(1+x),&x\in[-1,-B),
\end{cases}
\end{equation}
and the optimal control process $A^*\in\cA_1$ in the problem $V^*(x,1)$ is
given via~\eqref{eq:09122015a0.5} (with $a=1$) by the sequence of stopping
times
\begin{equation}\label{eq:09122015a9}
\begin{aligned}
&\tau_0^* = 0,\\ &\tau^*_{2n+1} = \inf\{t\ge \tau^*_{2n}: M_t \le
-B\},\\ &\tau^*_{2n+2} = \inf\{t\ge \tau^*_{2n+1}: M_t \ge B\}
\ \text{for}\ n\ge0.
\end{aligned}
\end{equation}

\noindent
(ii) If $c_1\ge (2\lambda + \alpha)^{-1}$, we have
$$
V^*(x,1)=\tilde V(x),\quad x\in[-1,1],
$$
and the optimal control process $A^*\in\cA_1$ is never to switch:
$A^*\equiv1$.
\end{theorem}

\begin{remark}
\label{rem:symmetry}
For the sake of comparison with~\eqref{eq:09122015a8}
and~\eqref{eq:09122015a9} let us again mention that $V^*(x,-1)=V^*(-x,1)$
and explicitly state that the optimal control $A^*\in\cA_{-1}$ in the
problem $V^*(x,-1)$ for the (interesting) case $c_1<(2\lambda +
\alpha)^{-1}$ is given via~\eqref{eq:09122015a0.5} (with $a=-1$) by the
sequence of stopping times
\begin{equation}\label{eq:09122017a2}
\begin{aligned}
&\tau_0^* = 0,\\ &\tau^*_{2n+1} = \inf\{t\ge \tau^*_{2n}: M_t \ge B\},\\
&\tau^*_{2n+2} = \inf\{t\ge \tau^*_{2n+1}: M_t \le -B\}\ \text{for}\ n\ge0.
\end{aligned}
\end{equation}
\end{remark}

\begin{remark}
A formal proof of \Cref{th:main} is performed in \Cref{sec:proof-main}. It
is based on a verification argument and does not explain how to come to the
statement of the theorem. The following discussion explains main ideas that
lead to the conditions appearing in it.

First, we reduce the optimization problem~\eqref{V} that contains the
process $\theta_t$, which is not $(\cF_t)$-adapted, to a problem containing
the posterior mean process $M_t$. It is not hard to see that the
distribution of $\theta_t$ conditionally on $\cF_t$ is given by the formula
$$
\Law(\theta_t|\cF_t)=\frac{1+M_t}2 \delta_1+\frac{1-M_t}2 \delta_{-1}.
$$
Since the controls $A$ are $(\cF_t)$-adapted and $\{-1,1\}$-valued, we get
$$
\PP_x(\theta_t\ne A_t|\cF_t)=\frac{1-A_t M_t}2.
$$
Taking intermediate conditioning with respect to $\cF_t$ in~\eqref{J} we
obtain for $x\in[-1,1]$, $A\in\cA$
\begin{multline}\label{J-repr}
J(x,A)=\EE_x\biggl[
\frac12\int_0^\infty e^{-\alpha t} (1-A_t M_t)\,dt
\\+\sum_{t\ge0} e^{-\alpha t} I(\Delta A_t\ne0)\Bigl(c_1 + \frac{c_2}{2}(1-A_tM_t)\Bigr)
\biggr].
\end{multline}
This provides a restatement of the optimisation problem~\eqref{V} in terms
of the process~$M$, which is introduced directly as a (unique strong)
solution to~\eqref{eq:04122015a3.5}.

It is intuitively clear that we can get a non-trivial solution only when the
switching costs are not too high, since otherwise it will never be optimal to
switch. This distinguishes cases (i) and (ii) of the theorem.

In order to solve the problem $V^*(x,1)$ in case (i), it is natural to guess
that the optimal strategy should be of the form~\eqref{eq:09122015a9} with
an appropriate threshold $B\in(0,1)$. Then we can expect that the value
function solves the inhomogeneous ODE~\eqref{ode} in $(-B,1)$, which is
known from the general optimal stopping theory (see, e.g.
\cite[Chapter~III]{PeskirShiryaev06}). The solutions to~\eqref{ode} that are
bounded in a left neighbourhood of $1$ are described
in~\eqref{eq:09122015a5}. Moreover, we are interested only in solutions with
$K>0$ in~\eqref{eq:09122015a5} because $\tilde V(x)$ turns out to be an
upper bound for $V^*(x,1)$. That can be seen from the equality
\begin{equation}\label{eq:09122015a1}
\tilde V(x)=\EE_x\left[
\frac12\int_0^\infty e^{-\alpha t} (1-M_t)\,dt\right]
\;\;(= J(x,1)),
\quad x\in[-1,1],
\end{equation}
which follows from direct computations.

Then it remains to find appropriate $K\in(0,\infty)$ and $B\in(0,1)$. If we
assume that the optimal strategy is of the form~\eqref{eq:09122015a9}, then
the value function $V^*(x,1)$ must be obtained by pasting like
in~\eqref{eq:09122015a8}. Natural conditions for determining $K$ and $B$ are
then continuous fit~\eqref{eq:09122015a6} and smooth
fit~\eqref{eq:09122015a7} at the point~$-B$.
\end{remark}

\subsection{Numerical solution}
Theorem~\ref{th:main} provides a complete solution of our problem. Now the
question is, given the numbers $\lambda,\mu,\alpha>0$, $c_1,c_2\ge0$ with
$c_1+c_2>0$ as inputs, how to determine the pair $(K,B)$ and the function
$\ffi$ numerically with a given precision. The problem is that $\ffi$ is not
known explicitly, but only characterised as a (unique up to a positive
multiplier) decreasing strictly positive solution of~\eqref{eq:hom_ode}.
This characterisation does not provide us with a Cauchy problem
for~\eqref{eq:hom_ode}, neither we can set up an appropriate boundary value
problem (recall that $\ffi(-1)=+\infty$). The following result about the
structure of $\ffi$ allows us, in particular, to write down a Cauchy problem
for~$\ffi$, which can be then efficiently solved numerically. This result
will be also used in the proof of Theorem~\ref{th:main}.

\begin{proposition}\label{th:phi}
For the function $\ffi$ from \Cref{pr:phi-properties} the following
statements hold true:
\begin{enumerate}[itemsep=0mm,label=(\alph*),leftmargin=*,topsep=\lineskip]
\item $\ffi$ is a strictly convex function on $(-1,1)$;
\item $\ffi'(1):=\lim\limits_{x\to 1}\ffi'(x)=-\frac{\alpha\ffi(1)}{2\lambda}<0$.
\end{enumerate}
\end{proposition}

\begin{figure}[p]
\centering
\includegraphics[width=0.327\textwidth]{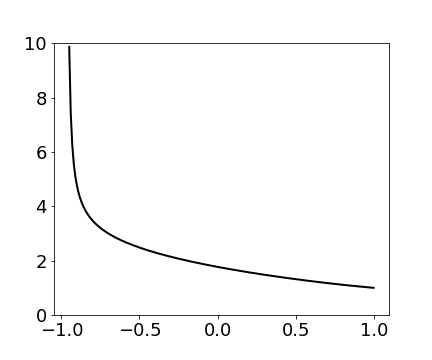}
\includegraphics[width=0.327\textwidth]{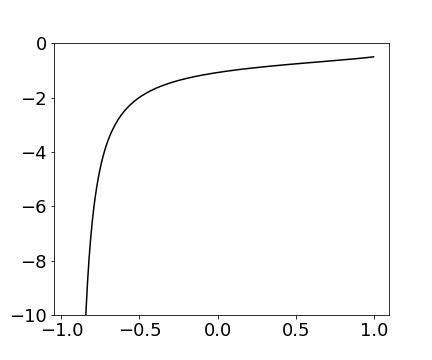}
\includegraphics[width=0.327\textwidth]{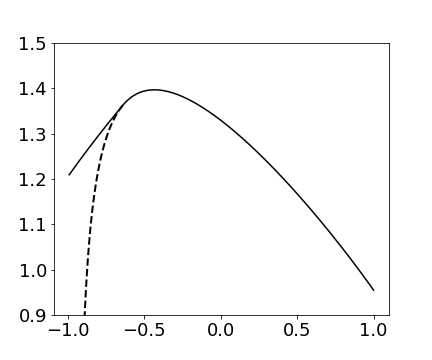}
\caption{\small Left: the function $\ffi$ with $\ffi(1)=1$ for
$\lambda=\alpha=1/4$, $\mu=1$.
Center: the derivative~$\ffi'$.
Right: the solid line is the value function $V^*(x,1)$; the dashed line is
$V(x)$, $x\in[-1,-B]$ (cf.~\eqref{eq:09122015a8}).
\label{fig:phi}}
\end{figure}

\begin{figure}[p]
\centering
\includegraphics[width=0.6\textwidth]{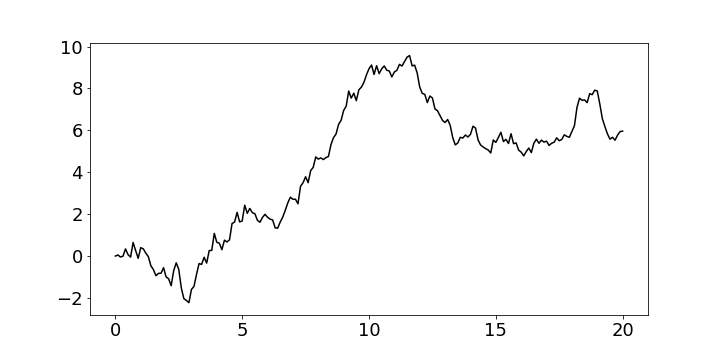}
\includegraphics[width=0.6\textwidth]{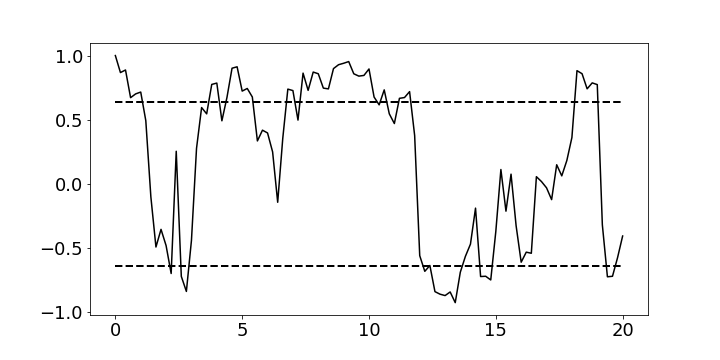}
\includegraphics[width=0.6\textwidth]{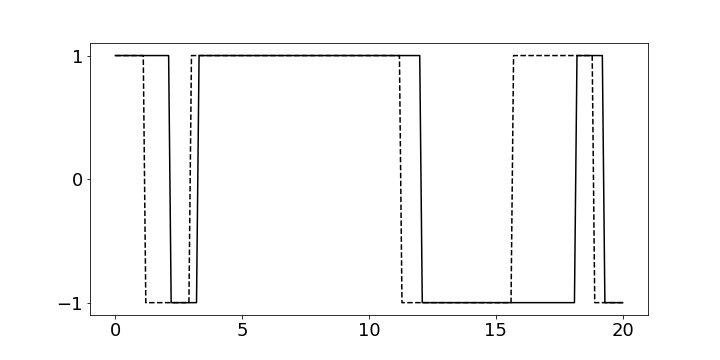}
\caption{\small Top: a simulated path of the process $X_t$.
Center: the posterior mean $M_t$ and the thresholds $\pm B$.
Bottom: the solid line is the control process $A_t$, the dashed line is the
unobservable process $\theta_t$.
\label{fig:simul}}
\end{figure}
In \Cref{fig:phi} we show the plot of the approximation of the function
$\ffi$ (normalised so that $\ffi(1)=1$) that corresponds to
$\lambda=\alpha=1/4$, $\mu=1$, and the plot of its derivative. The plots
were obtained by solving the Cauchy problem for~\eqref{eq:hom_ode} with
$\ffi(0.9999)=1$ and $\ffi'(0.9999)=-0.5$ (cf.~the second statement of
\Cref{th:phi}). We need to step a little to the left from $1$ in setting up
the Cauchy problem because \eqref{eq:hom_ode} is a regular ODE only strictly
inside $[-1,1]$ (the coefficient in front of the second derivative vanishes
at~$\pm1$).

Finally, being able to obtain good approximations of $\ffi$ we can
numerically find the values $(K,B)$ satisfying
\eqref{eq:09122015a6}--\eqref{eq:09122015a7} (see
\eqref{eq:16122015a1}--\eqref{eq:16122015a2} below for more detail). For
example, for the above values of $\lambda$, $\alpha$, $\mu$ and $c_1=1/4$,
$c_2=0$, we get $B=0.639$ and $K=0.378$ (this value of $K$ corresponds to
$\ffi$ normalised by $\ffi(1)=1$; recall that the value of $K$ depends on
the chosen version of~$\ffi$), which fully determines both the value
function $V^*(x,1)$ and the optimal strategy.

\Cref{fig:simul} shows a simulated example of the process $X_t$, the
posterior mean process $M_t$, and the optimal control process $A_t$ with the
same parameters.

\section{Proofs}\label{sec:proofs}
\subsection{Proofs of the auxiliary results}
\label{sec:proof-aux}
\begin{proof}[Proof of~\Cref{pr:phi-properties}]
While statements like (a)--(c) belong to common knowledge in the
Sturm-Liouville theory, some of their aspects depend on the ODE under
consideration. For example, for the ODE $f''-f=0$ on $(-1,1)$ (cf.\
with~\eqref{eq:hom_ode}), the functions $e^{-x}$ and $e^{-x+1}+e^{x-1}$ are
non-proportional strictly decreasing strictly positive solutions, that is,
the analogue of~(b) does not hold true. Therefore, even though the topic is
well-studied, we give a short treatment of claims (a)--(c) deducing them
from the general theory of one-dimensional diffusions. The analysis below
relies on boundary behavior of the process $M$, which is investigated in the
appendix.

Consider the diffusion $M$ driven by~\eqref{eq:04122015a3.5} inside $(-1,1)$
with a starting point $x\in(-1,1)$ under measure~$\PP_x$. Notice that the
first- and second-order derivative terms in~\eqref{eq:dif_op} constitute the
generator of~$M$, hence the relation between $M$ and ODE~\eqref{eq:hom_ode}.
As for claim~(a) of \Cref{pr:phi-properties}, the fact that
\eqref{eq:hom_ode} possesses a strictly decreasing (as well as a strictly
increasing) strictly positive solution is a well-known fact in the theory of
one-dimensional diffusions: e.g., see Proposition~50.3 in Chapter~V
in~\cite{RogersWilliams:00} or Section~II.1.10 in~\cite{BorodinSalminen:02}.
Moreover, there is a strictly decreasing strictly positive solution $\ffi$
of~\eqref{eq:hom_ode} with the property (analogue of formula~(50.5) in
Chapter~V of~\cite{RogersWilliams:00})
\begin{equation}\label{eq:10122017b1}
\ffi(x)=\ffi(y)\EE_x\left[e^{-\alpha T_y}\right],\quad \text{for all
}x>y\text{ in }(-1,1),
\end{equation}
where $\alpha>0$ is the one present in~\eqref{eq:dif_op} and
$$
T_y=\inf\{t\ge0:M_t=y\}
$$
(with the usual convention $\inf\emptyset:=\infty$).

Let $f\colon(-1,1)\to\bbR$ be a decreasing strictly positive solution
to~\eqref{eq:hom_ode}. In particular, $f\in C^2$. Fix a starting point
$x\in(-1,1)$ of $M$ and consider the process $(e^{-\alpha
t}f(M_t))_{t\ge0}$, which is well-defined because the boundaries $\pm1$ are
inaccessible for $M$, and is a $\PP_x$-local martingale by the It\^o
formula. Then, for any $y\in(-1,x)$, the process $X_t:=e^{-\alpha (t\wedge
T_y)}f(M_{t\wedge T_y})$, $t\ge0$, is a bounded $\PP_x$-martingale. The
identity $\EE_x X_0=\EE_x X_\infty$ now yields
$f(x)=f(y)\EE_x\left[e^{-\alpha T_y}\right]$. As $x>y$ in $(-1,1)$ are
arbitrary, we get that the function $f$ satisfies~\eqref{eq:10122017b1} (in
place of~$\ffi$). But~\eqref{eq:10122017b1} determines a function up to
proportionality (given a value of $f(c)$ for some fixed $c\in(-1,1)$, we
uniquely determine $f(x)$ for all $x\in(-1,1)$ via~\eqref{eq:10122017b1}).
This implies claim~(b).

Finally, for a fixed~$x$, let $y\to-1$ in~\eqref{eq:10122017b1}. Since the
boundary $-1$ is inaccessible for $M$, we have $\EE_x\left[e^{-\alpha
T_y}\right]\to0$ and, therefore, $\ffi(y)\to+\infty$ as $y\to-1$. Let now
$x\to1$ in~\eqref{eq:10122017b1} with a fixed~$y$. As $1$ is entrance
boundary, we have $\lim_{x\to1}\EE_x\left[e^{-\alpha T_y}\right]>0$, which
means that $\lim_{x\to1}\ffi(x)>0$. This proves claim~(c).
\end{proof}

\begin{proof}[Proof of~\Cref{th:phi}]
Given arbitrary strictly positive numbers $\lambda$, $\mu$ and $\alpha$,
take any decreasing strictly positive solution $\ffi$ to~\eqref{eq:hom_ode}.

\smallskip (a)  We have
\begin{equation}\label{eq:10122015b1}
\frac{\mu^2}2(1-x^2)^2\ffi''(x) =2\lambda x\ffi'(x)+\alpha\ffi(x), \quad
x\in(-1,1),
\end{equation}
hence $\ffi''>0$ on $(-1,0]$. Let us show that also $\ffi''\ge0$ on $[0,1)$.
Indeed, if there was a point $x_0\in[0,1)$ such that $\ffi''(x_0)< 0$, then
$\ffi'(x)$ would be locally decreasing at $x_0$, from which one would
conclude that $2\lambda x \ffi'(x) + \alpha \ffi(x) \le 2\lambda x
\ffi'(x_0) + \alpha \ffi(x_0) $ for any $x\in(x_0,1)$, which implies that
$\ffi''(x) \le \gamma/(1-x^2)^2$ for $x\in(x_0,1)$ with $\gamma =
(1-x_0^2)^2 \ffi''(x_0) < 0$. Integrating twice we would get that $\ffi(x)
\to -\infty$ as $x\to1$, which contradicts the positivity of~$\ffi$.

Thus, we have $\ffi''(x)\ge 0$ on $(-1,1)$, so $\ffi(x)$ is convex.
Moreover, it is actually \emph{strictly} convex on $(-1,1)$ because
otherwise there would be an interval, where $\ffi$ is affine, but there is
no affine function locally satisfying~\eqref{eq:hom_ode}.

\smallskip (b) Since $\ffi$ is decreasing and convex on $(-1,1)$, there is a
finite limit $\lim\limits_{x\to1}\ffi'(x)$ and
\begin{equation}\label{eq:10122015b2}
\lim_{x\to1}\ffi'(x)\le0.
\end{equation}
Then it follows from~\eqref{eq:10122015b1} that there is a finite limit
$\gamma:=\lim\limits_{x\to1}(1-x^2)^2\ffi''(x)$. Convexity of $\ffi$ yields
$\gamma\ge0$. If we assume that $\gamma>0$, we get
$\ffi''(x)\ge\frac12\gamma/(1-x^2)^2$ for $x\in[y,1)$ with some $y\in(0,1)$.
Integrating we obtain $\ffi'(x)\to+\infty$ as $x\to1$, which
contradicts~\eqref{eq:10122015b2}. Therefore,
$\lim\limits_{x\to1}(1-x^2)^2\ffi''(x)=0$, which concludes the proof.
\end{proof}

\subsection{Proof of the main theorem}
\label{sec:proof-main}
\newitem Consider the case $c_1<(2\lambda+\alpha)^{-1}$. It is convenient to
define $\beta=(2\lambda+\alpha)^{-1}$ and $\gamma =
\frac{c_1+{c_2}/{2}}{\beta+{c_2}/{2}}$. Let us prove that there is a pair
$(K,B)$ with $K>0$ and $B\in(\gamma,1)$ such that
\eqref{eq:09122015a6}--\eqref{eq:09122015a7} hold true, where we set
\begin{equation}\label{eq:10122015c0}
V(x)=\tilde V(x)-K\ffi(x),\quad x\in(-1,1]
\end{equation}
(recall that $\tilde V$ is given in~\eqref{eq:09122015a1} and that $\ffi$
satisfies claims (a)--(c) after~\eqref{eq:hom_ode}).

Let us observe that the system \eqref{eq:09122015a6}--\eqref{eq:09122015a7}
is equivalent to
\begin{equation}\label{eq:16122015a1}
K = h_1(B) = h_2(B),
\end{equation}
where the continuous functions $h_1$ and $h_2$ are defined on $(0,1)$ by the
formulas
\begin{equation}\label{eq:16122015a2}
h_1(x) = \frac{(\beta +{c_2}/{2})x-c_1-{c_2}/{2}}{\ffi(-x) - \ffi(x)},\quad
h_2(x) = \frac{\beta+c_2/2}{-\ffi'(-x) - \ffi'(x)}.
\end{equation}
Notice that not only $\ffi$ is strictly decreasing on $(-1,1)$, but also
$\ffi'<0$ everywhere on $(-1,1)$, so that we have strictly positive
quantities in the denominators. Indeed, if $\ffi'(x_0)=0$ at some point
$x_0\in(-1,1)$, then, by~\eqref{eq:hom_ode}, $\ffi''(x_0)>0$, hence $\ffi'$
would be strictly positive in a right neighbourhood of $x_0$. This would
contradict the fact that $\ffi$ is decreasing.

Since we are considering the case $c_1<\beta$, the function $h_1$ is
negative on $(0,\gamma)$ and positive on $(\gamma,1)$, while the function
$h_2$ is positive everywhere on $(0,1)$. Therefore, in order to prove that
the system \eqref{eq:09122015a6}--\eqref{eq:09122015a7} has a solution
$(K,B)\in(0,\infty)\times(\gamma,1)$, it is sufficient to show that
$h_1(x)/h_2(x) \ge 1$ for some $x\in(\gamma,1)$. Let $\eps>0$ be such that
$\gamma+\eps < 1$. Then for $x\in[\gamma + \eps,1)$ we have
$$
\frac{h_1(x)}{h_2(x)} = (x-\gamma)\,\frac{-\ffi'(-x) - \ffi'(x)}{\ffi(-x) -
\ffi(x)} \ge \eps\frac{-\ffi'(-x)}{\ffi(-x)},
$$
hence it will be enough to establish that
$\limsup_{x\to1}(-\ffi'(-x)/\ffi(-x)) =+\infty$, or, equivalently, with
$f(x):=\ffi'(x)/\ffi(x)$, that it holds
\begin{equation}\label{eq:10122015c2}
\liminf_{x\to-1}f(x)=-\infty.
\end{equation}
Dividing the homogeneous equation $L\ffi = 0$ by the function $\ffi$, we get
$$
\frac{\mu^2}{2}(1-x^2)^2 (f'(x) + f^2(x)) - 2\lambda xf(x) = \alpha,
$$
which implies
$$
f'(x) +f^2(x)=\frac{2(\alpha + 2\lambda xf(x))}{\mu^2(1-x^2)^2}.
$$
Observe that $f(x)\le0$, hence, if $x\le0$, then the right-hand side is
strictly positive and tends to $+\infty$ as $x\to-1$. If we suppose that $f$
is bounded in a right neighbourhood of~$-1$, then, for some $\eps>0$, we
have $f'(x)\ge\eps/(1-x^2)^2$ in a right neighbourhood of~$-1$, and
integrating we obtain a contradiction with the boundedness of~$f$, which
proves~\eqref{eq:10122015c2}. Thus, the system
\eqref{eq:09122015a6}--\eqref{eq:09122015a7} has a solution
$(K,B)\in(0,\infty)\times(\gamma,1)$. Its uniqueness will be proved later.

\newitem\label{it:10122017a1} Let us take any solution
$(K,B)\in(0,\infty)\times(\gamma,1)$ of the system
\eqref{eq:09122015a6}--\eqref{eq:09122015a7} and, with the function $V(x)$
given by~\eqref{eq:10122015c0}, define the ``candidate'' value function
$V(x,a)$, $x\in[-1,1]$, $a\in\{-1,1\}$, as the right-hand side
of~\eqref{eq:09122015a8}:
\begin{align}
V(x,1)&=\begin{cases}
V(x),&x\in[-B,1],\\
V(-x)+c_1+\dfrac{c_2}{2}(1+x),&x\in[-1,-B),
\end{cases}\label{eq:10122015c4}\\[1mm]
V(x,-1)&=V(-x,1),\quad x\in[-1,1].\label{eq:10122015c5}
\end{align}
We are going to prove that $V(x,a)$ defined in this way coincides with the
value function $V^*(x,a)$ defined in~\eqref{V}, where $J(x,A)$ admits
representation~\eqref{J-repr}.

Before we are able to do this, we need to establish the following auxiliary
facts:
\begin{enumerate}[itemsep=0mm,label=(F\arabic*),leftmargin=*,topsep=\lineskip]
\item  $V'(x,a)$ is bounded for $x\in[-1,1]$, where by
$V'(x,a)$ we will denote the derivative with respect to the first argument;
\item  $|\Delta V(x,a)|\le c_1+\frac{c_2}{2}(1+x)$ for
$x\in[-1,1]$, where $\Delta V(x,a):=V(x,a)-V(x,-a)$;
\item  $V(x,a)\in C^1$ in $x\in[-1,1]$, $V(x,a) \in C^2$ in
$x$ except at points $x=\pm B$, and $LV(x,a) \ge -\frac12 (1-ax)$ for
$x\in[-1,1]\setminus\{-B,B\}$.
\end{enumerate}

Due to the symmetry, it is enough to consider only $V(x,1)$. As
for~(F1), by~\eqref{eq:10122015c4} it is enough to establish boundedness of
$V'(x,1)$ in a left neighbourhood of~$1$. The latter follows
from~\eqref{eq:10122015c0} and the fact that $\ffi$ is decreasing and convex
(recall \Cref{th:phi}).

Let us prove that (F2). Denote for convenience $$ g(x):=\Delta V(x,1)
\;(\equiv V(x,1)-V(-x,1)).
$$
According to the definition of $V(x,a)$, we immediately have
$|g(x)|=c_1+\frac{c_2}{2}(1+x)$ for $|x|\ge B$. In particular,
$g(-B)=c_1+\frac{c_2}{2}(1-B)$, $g(B)=-c_1-\frac{c_2}{2}(1+B)$, so it is
sufficient to show that $g$ is decreasing on $[-B,B]$.

Considering the second derivative of $g$ in $(-B,B)$, $g''(x) =
-K(\ffi''(x)-\ffi''(-x))$, and recalling that (as follows
from~\eqref{eq:hom_ode})
\[
\ffi''(x)=\frac2{\mu^2}\, \frac{2\lambda
x\ffi'(x)+\alpha\ffi(x)}{(1-x^2)^2},
\]
hence, for $x<0$,
\[
\ffi''(x)\ge\frac2{\mu^2}\, \frac{\alpha\ffi(x)}{(1-x^2)^2}>\frac2{\mu^2}\,
\frac{\alpha\ffi(-x)}{(1-x^2)^2}\ge\ffi''(-x),
\]
we get $g''(x) < 0$ and $g'(x)$ is strictly decreasing for $x\in[-B,0)$. In
a similar way, $g'(x)$ is strictly increasing for $x\in (0,B]$. These two
facts combined with that $g'(-B) = g'(B) = 0$ (according
to~\eqref{eq:09122015a7}), mean that $g'(x)<0$ for $x\in(-B,B)$, so $g$ is
decreasing on $[-B,B]$. Thus, fact~(F2) is proved.

Finally, the differentiability properties of $V(x,1)$ follow from its
definition together with~\eqref{eq:09122015a7}, and, since $\tilde V(x)$
solves the inhomogeneous ODE, while $\ffi(x)$ solves the homogeneous ODE, we
have
\[
\begin{aligned}
&LV(x,1) = -\frac12(1-x) \text{ for }x\in(-B,1), \\
&LV(x,1) = -\frac12(1+x)-\alpha c_1 - \frac{c_2}{2}(2\lambda x + 1 + x)
\text{ for } x\in(-1,-B).
\end{aligned}
\]
Using that $c_1<\beta$ and $B>\gamma$ one can check that $LV(x,1) >
-\frac12(1-x)$ for $x\in(-1,-B)$, which completes the proof of~(F3).

\newitem\label{it:16122015a1} Now we continue with proving that $V(x,a)$
defined in \eqref{eq:10122015c4}--\eqref{eq:10122015c5} coincides with the
value function $V^*(x,a)$. Consider a control process $A_t$ starting in
$A_{0-}=a$ and let $(\tau_n)_{n=0}^\infty$ be the corresponding sequence of
stopping times (see \eqref{eq:09122015a0}--\eqref{eq:09122015a0.5}). By the
It\^o formula we have for any $n\ge 1$
\begin{multline*}
e^{-\alpha \tau_n}V(M_{\tau_n}, A_{\tau_n}) = V(M_0, a) + \sum_{1\le k\le
n}\biggl[ \int_{\tau_{k-1}}^{\tau_k} e^{-\alpha s} LV(M_s, A_s) ds \\+ \mu
\int_{\tau_{k-1}}^{\tau_k} e^{-\alpha s}(1-M_s^2) V'(M_s,A_s)\,d\ol W_s +
e^{-\alpha \tau_k} \Delta V(M_{\tau_k}, A_{\tau_k})\biggr].
\end{multline*}
Taking the expectation under $\PP_x$ of the both sides and letting $n\to
\infty$, so that $\tau_n\to\infty$, we obtain
\[
V(x,a) = -\E_x\sum_{k\ge 1}\biggl[ \int_{\tau_{k-1}}^{\tau_k} e^{-\alpha s}
LV(M_s, A_s) ds + e^{-\alpha \tau_k} \Delta V(M_{\tau_k},
A_{\tau_k})\biggr],
\]
where it was used that the stochastic integral is a uniformly integrable
martingale (since $V'(x,a)$ is bounded as proved in~(F1) above), so its
expectation is zero.

From (F2)--(F3) proved above, we get
\begin{multline}
V(x,a) \le \E_x \biggl[ \frac12 \int_0^\infty e^{-\alpha s} (1-A_sM_s)\,ds
\\+\sum_{s\ge 0}e^{-\alpha s} \I(\Delta A_s \neq 0) \Bigl(c_1 + \frac{c_2}{2}
(1-A_tM_t)\Bigr) \biggr].
\label{V-ineq}
\end{multline}
Taking the infimum over all $A\in \mathcal{A}_a$ we see that
\[
V(x,a) \le V^*(x,a), \quad x\in[-1,1].
\]
Furthermore, if we define the control process $A^*$ as
in~\eqref{eq:09122015a9} in the case $a=1$, respectively as
in~\eqref{eq:09122017a2} in the case $a=-1$, then we have
$LV(M_s,A^*_s)=-\frac12(1-A^*_sM_s)$ for all
$s\in(0,\infty)\setminus\{\tau_k:k\ge1\}$ and $\Delta V(M_{\tau_k},
A^*_{\tau_k}) = -(c_1 + \frac{c_2}{2}(1-A_s^* M_s))$ for all $k\ge1$, so
that \eqref{V-ineq} turns into equality, i.e., we get
\begin{multline*}
V(x,a) = \E_x \biggl[ \frac12 \int_0^\infty e^{-\alpha s} (1-A^*_sM_s)\,ds
\\+\sum_{s\ge 0}e^{-\alpha s} \I(\Delta A^*_s \neq 0)\Bigl(c_1 +
\frac{c_2}{2}(1-A_S^*M_s)\Bigr) \biggr].
\end{multline*}
Consequently, $V^*(x,a) = V(x,a)$.

\newitem The final step in considering the case $c_1<\beta$ is to prove
uniqueness of the pair $(K,B)\in(0,\infty)\times(\gamma,1)$ that satisfies
\eqref{eq:09122015a6}--\eqref{eq:09122015a7}, where the function $V(x)$,
$x\in[-1,1]$, is given by~\eqref{eq:10122015c0}. Suppose there are two such
pairs $(K_1,B_1)$ and $(K_2,B_2)$ and denote by $V_1(x,a)$, $V_2(x,a)$ the
respective functions defined by
\eqref{eq:10122015c4}--\eqref{eq:10122015c5}. Without loss of generality
assume $B_1\le B_2$.

According to the reasoning in part~\ref{it:16122015a1} of the proof, both
$V_1$ and $V_2$ coincide with the value function~$V^*$. In particular,
$V_1(x,1) = V_2(x,1)$ for $x\in[-B_1,1]$, which readily implies $K_1=K_2$.

Suppose $B_1<B_2$. Denote $V(x) = \tilde V(x) - K\ffi(x)$ with
$K=K_1\;(=K_2)$. It follows from~\eqref{eq:10122015c4} and the equality
$V_1(x,1)=V_2(x,1)$ that $V(x)=V(-x)+c_1+\frac{c_2}{2}(1+x)$ for
$x\in[-B_2,-B_1]$, i.e.,
\begin{equation}\label{eq:10122017a1}
g(x):=V(x)-V(-x)\equiv c_1+\frac{c_2}{2}(1+x)\quad\text{on }[-B_2,-B_1].
\end{equation}
However, by the treatment of~(F2) in part~\ref{it:10122017a1} of this proof,
$g''(x)<0$ for $x\in(-B_2,B_2)$, which contradicts~\eqref{eq:10122017a1}.
Therefore, $B_1=B_2$.

\newitem Finally, let us consider the case $c_1\ge\beta$. The proof of
part~(i) of Theorem~\ref{th:main} is completed by now, and we are going to
use it. To this end, we extend the previous notation to stress the
dependence of the involved functions and constants on~$c_1$. In particular,
for all $c_1\ge0$, $V^*(x,1)$ will now be denoted by $V^*_{c_1}(x,1)$; for
$c_1\in[0,\beta)$, the constants $K$ and $B$ determined by
\eqref{eq:09122015a5.5}--\eqref{eq:09122015a7} will be denoted by $K_{c_1}$
and~$B_{c_1}$ (we consider the same version of $\ffi$ for different values
of $c_1$, e.g., the version obtained by the normalisation $\ffi(1)=1$).
Also, for $c_1\in(0,\beta)$, we set
\begin{equation}\label{eq:10122017a2}
V_{c_1}(x)=\tilde V(x)-K_{c_1}\,\ffi(x),\quad x\in(-1,1]
\end{equation}
As $B_{c_1}\in(\gamma_{c_1},1)$ for
$c_1\in(0,\beta)$, we have $\lim_{c_1\to\beta}B_{c_1}=1$. Furthermore,
it follows from \eqref{eq:16122015a1}--\eqref{eq:16122015a2} that
$\lim_{c_1\to\beta}K_{c_1}=0$. Together with~\eqref{eq:09122015a8},
this yields
$$
\lim_{c_1\to\beta}V^*_{c_1}(x,1)=\tilde V(x), \quad x\in[-1,1].
$$
More precisely, for $x\in(-1,1]$, use \eqref{eq:09122015a8}, $B_{c_1}\to1$,
\eqref{eq:10122017a2}, and~$K_{c_1}\to0$, while, for $x=-1$, use
additionally that $\tilde V(-1)=\tilde V(1)+\beta$. Since, clearly,
$V^*_{c_1}(x,1)$ is increasing in~${c_1}$ (for any $x\in[-1,1]$), we obtain
$$
V^*_{c_1}(x,1)\ge\tilde V(x),\;\;x\in[-1,1],\quad \text{for all
}c_1\ge\beta.
$$
It remains to recall that, for any $c_1\ge\beta$, the strategy $A^*\in\cA_1$
with $A^*\equiv1$ produces the cost $\tilde V(x)$. This concludes the proof.

\appendix
\section*{Appendix: Boundary behaviour of the process~$M$}
In this appendix we show that the posterior mean process $M$ is a regular
diffusion in $(-1,1)$ with $\pm1$ being inaccessible entrance boundaries.
Such a characterization of the boundaries conceptually means that if a
solution to~\eqref{eq:04122015a3.5} is started at a point $x\in(-1,1)$, then
it never reaches the boundaries $\pm1$. If a solution is started at
$x=\pm1$, then it immediately enters the interval $(-1,1)$ and never leaves
it.

So, let us consider $M$ driven by~\eqref{eq:04122015a3.5} with a starting
point $x\in(-1,1)$. The regularity follows from that
$(1+|b(x)|)\sigma^{-2}(x)$ is a locally integrable function in $(-1,1)$,
where $b(x):=-2\lambda x$ is the drift coefficient and
$\sigma(x):=\mu(1-x^2)$ is the diffusion coefficient (see Section~5.5.C and,
in particular, conditions (ND)' and (LI)' in~\cite{KaratzasShreve:91}).

The scale function of $M$ is given by the formula
$$
p(x)=\int_0^x
\exp\left\{\frac{2\lambda}{\mu^2}\frac1{1-y^2}\right\}dy,
\quad x\in(-1,1).
$$
Since $p(\pm1)=\pm\infty$, the boundaries $\pm1$ are inaccessible
(see Proposition~5.5.22 in~\cite{KaratzasShreve:91}).

The speed measure of $M$ is
$$
m(dx)=\frac{2\,dx}{p'(x)\mu^2(1-x^2)^2}\,dx
\quad\text{on }(-1,1).
$$
Let us show that
\begin{equation}\label{eq:10122015b4}
\int_{(0,1)}p(y)\,m(dy)<\infty,
\end{equation}
which implies that $1$ is an entrance boundary for~$M$ (see Section~II.1.6
in~\cite{BorodinSalminen:02}), and, by symmetry, $-1$ is an entrance
boundary as well\footnote{Note that, in \cite{BorodinSalminen:02}, such a
boundary is called ``entrance-not-exit boundary'', while the term ``entrance
boundary'' in~\cite{BorodinSalminen:02} has a broader meaning. However, for
inaccessible boundaries both terminologies coincide.}.

In order to prove~\eqref{eq:10122015b4} it is enough to establish that there
is a finite limit
\[
\gamma:=\lim_{x\to1}\frac{p(x)}{p'(x)(1-x^2)^2}.
\]
We have
\[
p'(x)=\exp\left\{\frac{2\lambda}{\mu^2}\frac1{1-x^2}\right\},\qquad
p''(x)=p'(x)\frac{4\lambda}{\mu^2}\frac x{(1-x^2)^2}.
\]
Since the limit $\gamma$ is of the type $\frac\infty\infty$, by
l'H\^opital's rule, we find
$$
\gamma
=\lim_{x\to1}\frac{p'(x)}{p''(x)(1-x^2)^2
-2p'(x)(1-x^2)x}=\frac{\mu^2}{4\lambda}.
$$

\bibliographystyle{abbrvnat}
\bibliography{sequential_tracking}

\begin{thebibliography}{18}
\providecommand{\natexlab}[1]{#1}
\providecommand{\url}[1]{\texttt{#1}}
\expandafter\ifx\csname urlstyle\endcsname\relax
  \providecommand{\doi}[1]{doi: #1}\else
  \providecommand{\doi}{doi: \begingroup \urlstyle{rm}\Url}\fi

\bibitem[Bayraktar and Egami(2010)]{BayraktarEgami10}
E.~Bayraktar and M.~Egami.
\newblock On the one-dimensional optimal switching problem.
\newblock \emph{Mathematics of Operations Research}, 35\penalty0 (1):\penalty0
  140--159, 2010.

\bibitem[Bayraktar and Ludkovski(2009)]{BayraktarLudkovski09}
E.~Bayraktar and M.~Ludkovski.
\newblock Sequential tracking of a hidden {M}arkov chain using point process
  observations.
\newblock \emph{Stochastic Processes and their Applications}, 119\penalty0
  (6):\penalty0 1792--1822, 2009.

\bibitem[Borodin and Salminen(2002)]{BorodinSalminen:02}
A.~N. Borodin and P.~Salminen.
\newblock \emph{Handbook of {B}rownian Motion -- Facts and Formulae}.
\newblock Birkh\"auser Verlag, Basel, 2nd edition, 2002.

\bibitem[Brekke and {\O}ksendal(1994)]{BrekkeOksendal94}
K.~A. Brekke and B.~{\O}ksendal.
\newblock Optimal switching in an economic activity under uncertainty.
\newblock \emph{SIAM Journal on Control and Optimization}, 32\penalty0
  (4):\penalty0 1021--1036, 1994.

\bibitem[Cai et~al.(2017{\natexlab{a}})Cai, Rosenbaum, and
  Tankov]{CaiRosenbaumTankov17a}
J.~Cai, M.~Rosenbaum, and P.~Tankov.
\newblock Asymptotic lower bounds for optimal tracking: a linear programming
  approach.
\newblock \emph{Annals of Applied Probability}, 27\penalty0 (4):\penalty0
  2455--2514, 2017{\natexlab{a}}.

\bibitem[Cai et~al.(2017{\natexlab{b}})Cai, Rosenbaum, and
  Tankov]{CaiRosenbaumTankov17b}
J.~Cai, M.~Rosenbaum, and P.~Tankov.
\newblock Asymptotic optimal tracking: feedback strategies.
\newblock \emph{Stochastics}, 89\penalty0 (6-7):\penalty0 943--966,
  2017{\natexlab{b}}.

\bibitem[Duckworth and Zervos(2001)]{DuckworthZervos01}
K.~Duckworth and M.~Zervos.
\newblock A model for investment decisions with switching costs.
\newblock \emph{Annals of Applied probability}, 11\penalty0 (1):\penalty0
  239--260, 2001.

\bibitem[Gapeev(2015)]{Gapeev15}
P.~V. Gapeev.
\newblock Bayesian switching multiple disorder problems.
\newblock \emph{Mathematics of Operations Research}, 41\penalty0 (3):\penalty0
  1108--1124, 2015.

\bibitem[Karatzas and Shreve(1991)]{KaratzasShreve:91}
I.~Karatzas and S.~E. Shreve.
\newblock \emph{Brownian Motion and Stochastic Calculus}.
\newblock Springer-Verlag, New York, 2nd edition, 1991.

\bibitem[Liptser and Shiryaev(2001)]{LiptserShiryaev:01}
R.~S. Liptser and A.~N. Shiryaev.
\newblock \emph{Statistics of Random Processes {I}, {II}}.
\newblock Springer-Verlag, Berlin, 2001.

\bibitem[Ly~Vath and Pham(2007)]{LyVathPham07}
V.~Ly~Vath and H.~Pham.
\newblock Explicit solution to an optimal switching problem in the two-regime
  case.
\newblock \emph{SIAM Journal on Control and Optimization}, 46\penalty0
  (2):\penalty0 395--426, 2007.

\bibitem[Peskir and Shiryaev(2006)]{PeskirShiryaev06}
G.~Peskir and A.~Shiryaev.
\newblock \emph{Optimal Stopping and Free-Boundary Problems}.
\newblock Springer, 2006.

\bibitem[Pham(2009)]{Pham09}
H.~Pham.
\newblock \emph{Continuous-time Stochastic Control and Optimization with
  Financial Applications}.
\newblock Springer, 2009.

\bibitem[Poor and Hadjiliadis(2008)]{PoorHadjiliadis08}
H.~V. Poor and O.~Hadjiliadis.
\newblock \emph{Quickest Detection}.
\newblock Cambridge University Press, 2008.

\bibitem[Rogers and Williams(2000)]{RogersWilliams:00}
L.~C.~G. Rogers and D.~Williams.
\newblock \emph{Diffusions, {M}arkov Processes, and Martingales, vol.~2}.
\newblock Cambridge University Press, Cambridge, 2nd edition, 2000.

\bibitem[Shiryaev(2010)]{Shiryaev10}
A.~N. Shiryaev.
\newblock Quickest detection problems: fifty years later.
\newblock \emph{Sequential Analysis}, 29\penalty0 (4):\penalty0 445--385, 2010.

\bibitem[Shiryaev(2019)]{Shiryaev19}
A.~N. Shiryaev.
\newblock \emph{Stochastic Disorder Problems}.
\newblock Springer, 2019.

\bibitem[Tartakovsky et~al.(2014)Tartakovsky, Nikiforov, and
  Basseville]{TartakovskyNikiforov14}
A.~Tartakovsky, I.~Nikiforov, and M.~Basseville.
\newblock \emph{Sequential Analysis: Hypothesis Testing and Changepoint
  Detection}.
\newblock CRC Press, 2014.

\end{thebibliography}
\end{document}